\documentclass[12pt,leqno,a4paper]{article}
\usepackage{amsmath,amsthm,enumerate}
\usepackage[latin1]{inputenc}
\usepackage[T1]{fontenc}
\usepackage[english]{babel}

\newtheorem{theorem}{Theorem}[section]

\newtheorem{corollary}[theorem]{Corollary}

\newtheorem{remark}[theorem]{Remark}


\newcommand*\proofnamestyle{\itshape}

\DeclareMathOperator{\tr}{Tr}

\begin{document}

    \title{Jensen's inequality for conditional expectations}
    
    \date{March 22 2005}
    
    \author{Frank Hansen}
    \maketitle

    \begin{abstract} We study conditional expectations generated by an
    abelian $ C^* $-subalgebra in the centralizer of a positive functional. We formulate and prove
    Jensen's inequality for functions of several variables with respect to this type of conditional
    expectations, and we obtain as a corollary Jensen's inequality for expectation values.
    \end{abstract}

    \section{Preliminaries} 
    
    An $ n $-tuple $ \underline{x}=(x_1,\dots,x_n) $ of elements
    in a $ C^* $-algebra $ {\mathcal A} $ is said to be abelian if the elements 
    $ x_1,\dots,x_n $ are mutually commuting. We say that an abelian $ n $-tuple $ \underline{x} $
    of self-adjoint elements is in the domain of  a real continuous function $ f $ of $ n $ variables
    defined on a cube of real intervals $ \underline{I}=I_1\times\cdots\times I_n $ 
    if the spectrum $ \sigma(x_i) $ of $ x_i $ is contained in $ I_i $ for each $ i=1,\dots,n. $   
    In this situation $ f(\underline{x}) $ is naturally defined as an element in $ {\mathcal A} $
    in the following way. We may assume that $ {\mathcal A} $ is realized as operators on a Hilbert space
    and let
    \[
    x_i=\int\lambda\, dE_i(\lambda)\qquad i=1,\dots,n
    \]
    denote the spectral resolutions of the operators $ x_1,\dots,x_n. $
    Since the $ n $-tuple $ \underline{x}=(x_1,\dots,x_n) $ is abelian, the spectral
    measures $ E_1,\dots,E_n $ are mutually commuting. We may thus set
    \[
    E(S_1\times\cdots\times S_n)=E_1(S_1)\cdots E_n(S_n) 
    \]
    for Borel sets $ S_1,\dots,S_n $ in $ {\mathbf R} $ and extend $ E $ to a spectral measure 
    on $ {\mathbf R}^n $ with support in $ \underline{I}. $ Setting
    \[
    f(\underline{x})=\int f(\lambda_1,\dots,\lambda_n)\, dE(\lambda_1,\dots,\lambda_n)
    \]
    and since $ f $ is continuous,
    we finally realize that $ f(\underline{x}) $ is an element in $ {\mathcal A}. $
    
    \section{Conditional expectations}

    Let $ {\mathcal C} $ be a separable abelian $ C^* $-subalgebra of a $ C^* $-algebra $ {\mathcal A}, $ and let
    $ \varphi $ be a positive functional on $ {\mathcal A} $ such that $ {\mathcal C} $ is contained in the
    centralizer
    \[
    {\mathcal A}^{\varphi}=\{y\in{\mathcal A}\mid \varphi(xy)=\varphi(yx)\quad\forall x\in{\mathcal A}\}.
    \] 
    The subalgebra is of the form $ {\mathcal C}=C_0(S) $ for some locally compact metric space $ S. $
    
    \begin{theorem}
    There exists a positive linear mapping
    \begin{gather}\label{conditional expectation}
    \Phi\colon M({\mathcal A})\to L^\infty(S,\mu_\varphi)
    \end{gather}
    on the multiplier algebra $ M({\mathcal A}) $ such that
    \[
    \Phi(xy)=\Phi(yx)=\Phi(x)y\qquad x\in M({\mathcal A}),\, y\in {\mathcal C}
    \]
    almost everywhere, and a finite Radon measure $ \mu_\varphi $ on $ S $ such that
    \begin{gather}
    \int_S z(s)\Phi(x)(s)\, d\mu_\varphi(s)=\varphi(zx)\qquad z\in {\mathcal C},\, x\in M({\mathcal A}).
    \end{gather}
    \end{theorem}
    
    \begin{proof} By the Riesz representation theorem there is a finite
    Radon measure $ \mu_\varphi $ on $ S $ such that
    \[
    \varphi(y)=\int_S y(s)\,d\mu_\varphi(s)\qquad y\in{\mathcal C}=C_0(S).
    \]    
    For each positive element $ x $ in the multiplier algebra
    $ M(\mathcal A) $ we have
    \[
    0\le\varphi(yx)=\varphi(y^{1/2}xy^{1/2})\le\|x\|\varphi(y)\qquad y\in{\mathcal C}_+.
    \]
    The functional $ y\to\varphi(yx) $ on $ {\mathcal C} $ consequently defines a Radon measure on $ S $ which is
    dominated by a multiple of $ \mu_\varphi, $ and it is therefore given by a unique element $ \Phi(x) $ in
    $ L^\infty(S,\mu_\varphi). $ By linearization this defines a positive linear mapping defined on the
    multiplier algebra
    \begin{gather}\label{conditional expectation}
    \Phi\colon M({\mathcal A})\to L^\infty(S,\mu_\varphi)
    \end{gather}
    such that
    \[
    \int_S z(s)\Phi(x)(s)\, d\mu_\varphi(s)=\varphi(zx)\qquad z\in {\mathcal C},\, x\in M({\mathcal A}).
    \]
    Furthermore, since
    \[
    \int_S z(s)\Phi(yx)(s)\, d\mu_\varphi(s)=\varphi(zyx)=\int_S z(s)y(s)\Phi(x)(s)\, d\mu_\varphi(s)
    \]
    for $ x\in M({\mathcal A}) $ and $ z,y\in {\mathcal C} $ we derive $ \Phi(yx)=y\Phi(x)=\Phi(x)y $
    almost everywhere.
    Since $ {\mathcal C} $ is contained in the centralizer $ {\mathcal A}^\varphi $ and thus
    $ \varphi(zxy)=\varphi(yzx), $ we similarly obtain $ \Phi(xy)=\Phi(x)y $ almost everywhere. \hfill
    \end{proof}\vskip 2ex

    Note that $ \Phi(z)(s)=z(s) $ almost everywhere in $ S $ for each $ z\in {\mathcal C}, $
    cf. \cite{kn:lieb:2002, kn:hansen:2003:2, kn:hansen:2003:3}. 
    With a slight abuse of language we call $ \Phi $ a conditional expectation
    even though its range is not a subalgebra of $ M({\mathcal A}). $

    \section{Jensen's inequality}\label{section: expectation values}
    
    Following the notation in \cite{kn:hansen:2003:3} we consider a separable $ C^* $-algebra
    $ {\mathcal A} $ of operators on a (separable) Hilbert space $ H, $ and a field
    $ (a_t)_{t\in T} $ of operators in the multiplier algebra 
    \[
    M(\mathcal A)=\{a\in B(H)\mid a{\mathcal A} + {\mathcal A}a \subseteq {\mathcal A}\}
    \]
    defined on a locally compact metric space $ T $ equipped with a Radon measure $ \nu $.
    We say that the field $ (a_t)_{t\in T} $ is weak*-measurable if the function $ t\to\varphi (a_t) $ is 
    $ \nu $-measurable on $ T $ for each
    $ \varphi \in {\mathcal A}^*; $ and we say that the field is continuous if the function  $ t\to a_t $
    is continuous \cite{kn:hansen:2003:2}. 
    
    As noted in \cite{kn:hansen:2003:3} the field $ (a_t)_{t\in T} $ is weak*-measurable,
    if and only if for each vector $ \xi\in H $ the function $ t\to a_t\xi $ is weakly (equivalently
    strongly) measurable. In particular, the composed field $ (a_t^*b_t)_{t\in T} $ is weak*-measurable if
    both $ (a_t)_{t\in T} $ and $ (b_t)_{t\in T} $ are weak*-measurable fields.
    
    If for a weak*-measurable field $ (a_t)_{t\in T} $ the function $ t\to|\varphi (a_t)| $ is integrable for 
    every state $ \varphi\in S({\mathcal A}) $ and the integrals
    \[
    \int_T |\varphi(a_t)|\,d\nu(t)\le K\qquad\forall\varphi\in S({\mathcal A})
    \]
    are uniformly bounded by some constant $ K, $ then there is a unique element 
    (a $ C^*$-integral in Pedersen's terminology \cite[2.5.15]{kn:pedersen:1989}) in the multiplier algebra 
    $ M(\mathcal A), $ designated by 
    \[
    \int_T a_t\,d\nu(t),
    \]
    such that  
    \[
    \varphi \left(\int_T a_t\, d\nu (t)\right) = 
    \int_T \varphi (a_t)\, d\nu (t)  \qquad \forall\varphi \in {\mathcal A}^*.
    \]
    We say in this case that the field $(a_t)_{t\in T}$ is integrable. Finally we say that a field $ (a_t)_{t\in T} $ 
    is a unital column field \cite{kn:araki:2000,kn:hansen:2003:2, kn:hansen:2003:3}, if
    it is weak*-measurable and
    \[
    \int_T a_t^*a_t\, d\nu(t)=1.
    \]
    We note that a $ C^* $-subalgebra of a separable $ C^* $-algebra is automatically separable.    
    
     \begin{theorem}\label{theorem: Jensen's inequality for conditional expectations}
     Let $ {\mathcal C} $ be an abelian $ C^* $-subalgebra of a separable $ C^* $-algebra $ {\mathcal A}, $
     $ \varphi $ be a positive functional on $ {\mathcal A} $ such that $ {\mathcal C} $ is contained in the
    centralizer $ {\mathcal A}^\varphi $ and let 
    \[
    \Phi\colon M({\mathcal A})\to L^\infty(S,\mu_\varphi) 
    \]
    be the conditional expectation defined in (\ref{conditional expectation}). Let furthermore
    $ f:\underline{I}\to {\mathbf R} $ be a continuous convex function of $ n $ variables defined on a cube,
    and let $ t\to a_t\in M({\mathcal A}) $ be a unital column field on a locally compact Hausdorff space $ T $
    with a Radon
    measure $ \nu. $ If $ t\to\underline{x}_t $ is an essentially bounded, weak* measurable field on $ T $ of abelian 
    $ n $-tuples of self-adjoint elements in $ {\mathcal A} $ in the domain of $ f, $ then
    \begin{gather}\label{Jensen's inequality for conditional expectations}
    f(\Phi(y_1),\dots,\Phi(y_n))\le \Phi\left(\int_T a_t^*f(\underline{x}_t) a_t\,d\nu(t)\right)
    \end{gather}
    almost everywhere,
    where the $ n $-tuple $ \underline{y} $ in $ M({\mathcal A}) $ is defined by setting
    \[
    \underline{y}=(y_1,\dots,y_n)=\int_T a_t^* \underline{x}_t a_t\,d\nu(t).    
    \]
    \end{theorem}
    
    \begin{proof}  The subalgebra $ {\mathcal C} $ is as noted above of the form $ {\mathcal C}=C_0(S) $
     for some locally compact metric space $ S, $ and since the $ C^*$-algebra $C_0(\underline{I})$ is 
separable we may for almost every $s$ in $S$ define a 
Radon measure $\mu_s$ on $\underline{I}$ by setting
\[
\mu_s(g)=\int_{\underline{I}} g(\underline{\lambda})\, d\mu_s (\underline{\lambda}) = \Phi\left( \int_T a^*_t
g(\underline{x}_t)a_t\,d\mu(t)\right)(s) \qquad g\in C_0(\underline{I}).
\]
Since
\[
\mu_s(1)=\Phi\left(\int_T a^*_ta_t\,d\mu(t)\right)=\Phi(1)=1
\]
we observe that $ \mu_s $ is a probability measure. If we put $ g_i(\underline{\lambda})= \lambda_i $ then 
\[
\int_{\underline{I}} g_i(\underline{\lambda})\,d\mu_s(\underline{\lambda}) = \Phi\left(\int_T a^*_t x_{it} a_t\,
d\mu(t)\right)(s) = \Phi(y_i) (s)
\]
for $ i=1,\dots, n $ and since $ f $ is convex we obtain
\[
\begin{array}{l}
f(\Phi(y_1)(s),\dots, \Phi(y_n)(s))
=\displaystyle f\left(\int_{\underline{I}} g_1(\underline{\lambda})\,d\mu_s(\underline{\lambda}), \dots,
\int_{\underline{I}} g_n(\underline{\lambda})\, d\mu_s(\underline{\lambda})\right)\\[3ex]
\displaystyle\le\int_{\underline{I}} f\left(g_1(\underline{\lambda}), \dots , g_n(\underline{\lambda})\right)
\, d\mu_s(\underline{\lambda})=\int_{\underline{I}} f(\underline{\lambda})\,d\mu_s(\underline{\lambda})\\[3ex]
=\displaystyle\Phi\left(\int_T a^*_t f(\underline{x}_t)a_t\,d\mu(t)\right)(s)
\end{array}
\]
for almost all $ s $ in $ S. $\hfill\end{proof}\vskip 2ex
    
    The following corollary is known as ''Jensen's inequality for expectation values''. It was formulated
    (for continuous fields) in the reference \cite{kn:hansen:2005:4}, where a more direct proof is given.

    \begin{corollary}\label{Jensen's inequality for expectation values}
    Let $ f:\underline{I}\to {\mathbf R} $ be a continuous convex function of $ n $ variables defined on a cube,
    and let $ t\to a_t\in B(H) $ be a unital column field on a locally compact Hausdorff space $ T $ with a Radon
    measure $ \nu. $ If $ t\to\underline{x}_t $ is a bounded weak*-measurable field on $ T $ of abelian 
    $ n $-tuples of self-adjoint operators on $ H $ in the domain of $ f, $ then
    \begin{gather}\label{vector inequality}
    f\bigl((y_1\xi\mid\xi),\dots,(y_n\xi\mid\xi)\bigr)\le\left( 
    \int_T a_t^* f(\underline{x}_t) a_t\, d\nu(t)\xi\mid\xi\right) 
    \end{gather}
    for any unit vector $ \xi\in H, $ where the $ n $-tuple $ \underline{y} $ is defined by setting
    \[
    \underline{y}=(y_1,\dots,y_n)=\int_T a_t^* \underline{x}_t a_t\,d\nu(t).
    \]
    \end{corollary}
    
    \begin{proof}  
    The statement follows from Theorem \ref{theorem: Jensen's inequality for conditional expectations}
    by choosing $ \varphi $ as the trace and letting $ {\mathcal C} $ be the $ C^* $-algebra generated by the
    orthogonal projection $ P $ on the vector $ \xi. $  Then $ {\mathcal C}=C_0(S) $ where $ S=\{0,1\}, $ and an
    element $ z\in{\mathcal C} $ has the representation
    \[
    z=z(0)P+z(1)(1-P).
    \]    
    The measure $ d\mu_\varphi $ gives unit weight in each of the two points, and the conditional expectation 
    $ \Phi $ is given by
    \[
    \Phi(x)(s)=\left\{\begin{array}{ll}
                      (x\xi\mid\xi)   &s=0\\[1ex]
                      \tr(x-Px)\quad &s=1.
                      \end{array}\right.
    \]
    Indeed,
    \[
    \begin{array}{rl}
    \varphi(zx)&=\tr\Bigl(\bigl(z(0)P+z(1)(1-P)\bigr)x\Bigr)\\[2ex]
    &=z(0)\Phi(x)(0)+z(1)\Phi(x)(0)\\[2ex]
    &=\displaystyle\int_S z(s)\Phi(x)(s)\,ds
    \end{array}    
    \]
    as required. The statement follows by evaluating the functions appearing on each side of
    the inequality (\ref{Jensen's inequality for conditional expectations}) in the point $ s=0. $
    \end{proof} 
    
     \begin{remark}\label{remark}
    If we choose $ \nu $ as a probability measure on $ T, $ then the trivial field $ a_t=1 $
    for $ t\in T $ is unital and (\ref{vector inequality}) takes the form
    \[
    f\left(\left(\int_T x_{1t} \,d\nu(t)\xi\mid\xi\right),
    \dots, \left(\int_T x_{nt} \,d\nu(t)\xi\mid\xi\right)\right)\le
    \left(\int_T f(\underline{x}_t)\,d\nu(t) \xi\mid\xi\right)
    \]
    for bounded weak*-measurable fields of abelian $ n $-tuples $ \underline{x}_t=(x_{1t},\dots,x_{nt}) $
    of self-adjoint operators in the domain of $ f $  and unit vectors $ \xi. $
    By choosing $ \nu $ as an atomic measure with one atom we get a version 
    \begin{gather}\label{mond and pecaric}
    f\bigl(\left(x_1\xi\mid\xi\right),\dots, \left(x_n\xi\mid\xi\right)\bigr)\le
    \bigl(f(\underline{x})\xi\mid\xi\bigr)
    \end{gather}
    of the Jensen inequality by Mond and Pe{\v{c}}ari{\'{c}} \cite{kn:mond:1993:1}. By further considering a direct
    sum 
    \[
    \xi=\displaystyle\bigoplus_{j=1}^m \xi_j\quad\mbox{and}\quad
    x=(x_1,\dots,x_n)=\bigoplus_{j=1}^m (x_{1j},\dots,x_{nj})
    \]
    we obtain the familiar version
    \[
    f\left(\sum_{j=1}^m (x_{1j}\xi_j\mid \xi_j),\dots,\sum_{j=1}^m (x_{nj}\xi_j\mid \xi_j)\right)
    \le\sum_{j=1}^m \Bigl(f(x_{1j},\dots,x_{nj})\xi_j\mid\xi_j\Bigr)
    \]
    valid for abelian $ n $-tuples $ (x_{1j},\dots,x_{nj}), $ $ j=1,\dots,m $ of self-adjoint operators
    in the domain of $ f $ and vectors $ \xi_1,\dots,\xi_m $ with $ \|\xi_1\|^2+\cdots+\|\xi_m\|^2=1. $
    \end{remark}

      \vfill

      {\small\noindent Frank Hansen: Institute of Economics, University
       of Copenhagen, Studiestraede 6, DK-1455 Copenhagen K, Denmark.}

      \end{document}